\documentclass[11pt,leqno]{article}
\usepackage{amsmath,amssymb,amscd,latexsym}
\usepackage{calc}
\usepackage[all]{xy}

\newcommand{\C}{{\mathbb{C}}}
\newcommand{\F}{{\mathbb{F}}}
\newcommand{\oF}{\overline{\F}}
\newcommand{\Ge}{\mathbb{G}}

\newcommand{\Q}{{\mathbb{Q}}}
\newcommand{\oQ}{\overline{\Q}}
\newcommand{\Z}{{\mathbb{Z}}}
\newcommand{\oZ}{\overline{\Z}}

\newcommand{\bV}{\mathbf{V}}

\newcommand{\tC}{\tilde{C}}

\newcommand{\abb}{\mathrm{ab}}

\newcommand{\car}{\mathrm{char}}
\newcommand{\cont}{\mathrm{cont}}
\newcommand{\et}{\mathrm{\acute{e}t}}

\newcommand{\fin}{\mathrm{fin}}

\newcommand{\Mor}{\mathrm{Mor}}
\newcommand{\red}{\mathrm{red}}
\newcommand{\rk}{\mathrm{rk}\,}
\newcommand{\spec}{\mathrm{spec}\,}

\newcommand{\Aut}{\mathrm{Aut}}

\newcommand{\Ext}{\mathrm{Ext}}
\newcommand{\Fr}{\mathrm{Fr}}

\newcommand{\Gal}{\mathrm{Gal}}
\newcommand{\GL}{\mathrm{GL}\,}
\newcommand{\bGL}{\mathbf{GL}}
\newcommand{\Hom}{\mathrm{Hom}}

\newcommand{\Imm}{\mathrm{Im}\,}
\newcommand{\Pic}{\mathrm{Pic}}

\newcommand{\uVec}{\mathrm{Vec}}
\newcommand{\tors}{\mathrm{tors}}
\newcommand{\Ah}{{\mathcal A}}

\newcommand{\Ch}{{\mathcal C}}
\newcommand{\Dh}{{\mathcal D}}
\newcommand{\Eh}{{\mathcal E}}
\newcommand{\tEh}{\tilde{\Eh}}
\newcommand{\Fh}{{\mathcal F}}

\newcommand{\Oh}{{\mathcal O}}

\newcommand{\Th}{{\mathcal T}}
\newcommand{\Yh}{\mathcal{Y}}
\newcommand{\Zh}{\mathcal{Z}}

\newcommand{\emm}{{\mathfrak{m}}}
\newcommand{\eo}{\mathfrak{o}}

\newcommand{\eB}{\mathfrak{B}}
\newcommand{\eT}{\mathfrak{T}}
\newcommand{\eX}{{\mathfrak X}}

\newcommand{\eZ}{\mathfrak{Z}}

\newcommand{\tX}{\tilde{X}}
\newcommand{\silo}{\stackrel{\sim}{\longrightarrow}}

\newcommand{\verk}{\raisebox{0.03cm}{\mbox{\scriptsize $\,\circ\,$}}}
\newtheorem{theorem}{Theorem}
\newtheorem{lemma}[theorem]{Lemma}
\newtheorem{prop}[theorem]{Proposition}
\newtheorem{defn}[theorem]{Definition}
\newtheorem{cor}[theorem]{Corollary}

\newenvironment{proof}{\noindent {\bf Proof}}{\mbox{}\hspace*{\fill}$\Box$}
\newenvironment{proofob}{\noindent {\bf Proof}}{}
\newcounter{hours} \newcounter{minutes}

\begin{document}
\title{On Tannaka duality for vector bundles on $p$-adic curves}
\author{Christopher Deninger \and Annette Werner}
\date{\it  Dedicated to Jacob Murre}
\maketitle

\section{Introduction}
\label{sec:1}

In our paper \cite{De-We2} we have introduced a certain category $\eB^{ps}$ of degree zero bundles with ``potentially strongly semistable reduction'' on a $p$-adic curve. For these bundles it was possible to establish a partial $p$-adic analogue of the classical Narasimhan--Seshadri theory for semistable vector bundles of degree zero on compact Riemann surfaces. One of the main open questions of \cite{De-We2} was, whether our category was abelian. The first main result of the present note, Corollary \ref{t3_5}, asserts that this is indeed the case. It follows that $\eB^{ps}$ is even a neutral Tannakian category. In the second main result, theorem \ref{t4_6}, we calculate the group of connected components of the Tannaka dual group of $\eB^{ps}$. This uses a result of A. Weil characterizing vector bundles that become trivial in a finite \'etale covering as the ones satisfying a ``polynomial equation'' over the integers.

Bsides, in section \ref{sec:3} we give a short review of \cite{De-We2}, and in section \ref{sec:2} we discuss the ``strongly semistable reduction'' condition.

We would like to draw the reader's attention to the paper of Faltings \cite{Fa} on non-abelian $p$-adic Hodge theory which generalizes several results of \cite{De-We2} to Higgs bundles.

It is a pleasure for us to thank Uwe Jannsen for a helpful discussion.

\section{Vector bundles in characteristic $p$}
\label{sec:2}

Throughout this paper, we call a purely one-dimensional separated scheme of finite type over a field $k$ a curve over $k$.\\
Let $C$ be a connected smooth, projective curve over $k$. For a vector bundle $E$ on $C$ we denote by $\mu (E) = \frac{\deg (E)}{\rk (E)}$ the slope of $E$. Then $E$ is called semistable (respectively stable), if for all proper non-trivial subbundles $F$ of $E$ the inequality $\mu (F) \le \mu (E)$ (respectively $\mu (F) < \mu (E)$) holds. 

\begin{lemma}
  \label{t2_1}
If $\pi : C' \to C$ is a finite separable morphism of connected smooth projective $k$-curves, then semistability of $E$ is equivalent to semistability of $\pi^* E$.
\end{lemma}

\begin{proof}
  See \cite{Gie2}, 1.1.
\end{proof}

If $\car (k) = 0$, then by lemma \ref{t2_1} any finite morphism of smooth connected projective curves preserves semistability. However in the case $\car (k) = p$, there exist vector bundles which are destabilized by the Frobenius map, see \cite{Gie1}, Theorem 1. Assume that $\car (k) = p$, and let $F : C \to C$ be the absolute Frobenius morphism, defined by the $p$-power map on the structure sheaf. 

\begin{defn}
  \label{t2_2}
A vector bundle $E$ on $C$ is called strongly semistable of degree zero if $\deg (E) = 0$ and if $F^{n*} E$ is semistable on $C$ for all $n \ge 0$. 
\end{defn}

Now we also consider non-smooth curves over $k$. Let $Z$ be a proper curve over $k$. By $C_1 , \ldots , C_r$ we denote the irreducible components of $Z$ endowed with their reduced induced structures. Let $\tC_i$ be the normalization of $C_i$, and write $\alpha_i : \tC_i \to C_i \to Z$ for the canonical map. Note that the curves $\tC_i$ are smooth irreducible and projective over $k$.

\begin{defn}
  \label{t2_3}
A vector bundle $E$ on the proper $k$-curve $Z$ is called strongly semistable of degree zero, if all $\alpha^*_i E$ are strongly semistable of degree zero.
\end{defn}

A alternative characterization of this propery is given by the following result.

\begin{prop}
  \label{t2_4}
A vector bundle $E$ on $Z$ is strongly semistable of degree zero if and only if for any $k$-morphism $\pi : C \to Z$, where $C$ is a smooth connected projective curve over $k$, the pullback $\pi^* E$ is semistable of degree zero on $C$.
\end{prop}

Note that in \cite{De-Mi}, (2.34) bundles with this property are called semistable of degree zero.

\begin{proof}
  Let $X$ be a scheme over $k$. The absolute Frobenius $F$ sits in a commutative diagram
\[
\xymatrix{
X \ar[r]^F \ar[d] & X \ar[d] \\
\spec k \ar[r]^F & \spec k
}
\]
If we denote for all $r \ge 1$ by $X^{(r)}$ the scheme $X$ together with the structure map $X \to \spec k \xrightarrow{F^r} \spec k$, then $F^r : X^{(r)} \to X$ is a morphism over $\spec k$.

Assume that the curve $Z$ has the property in the claim. Applying it to the smooth projective curves $\tC_i$ and the $k$-morphisms
\[
\tC^{(r)}_i \xrightarrow{F^r} \tC_i \xrightarrow{\alpha_i} Z \; ,
\]
we find that all $\alpha^*_i E$ are strongly semistable of degree zero, i.e. that $E$ is strongly semistable of degree zero in the sense of definition \ref{t2_3}.

Conversely, assume that all $\alpha^*_i E$ are strongly semistable of degree zero. Let $\pi : C \to Z$ be a $k$-morphism from a smooth connected projective curve $C$ to $Z$. Then $\pi$ factors through one of the $C_i$. If $\pi$ is constant, then $\pi^* E$ is trivial, hence semistable of degree zero. Hence we can assume that $\pi (C) = C_i$. Since $C$ is smooth, $\pi$ also factors through the normalization $\tC_i$, i.e. there is a morphism $\pi_i : C \to \tC_i$ satisfying $\alpha_i \verk \pi_i = \pi$. Since $\pi$ is dominant, it is finite and hence $\pi_i$ is the composition of a separable map and a power of Frobenius, see e.g. \cite{Ha}, IV, 2.5. Hence there exists a smooth projective curve $D$ over $k$ and a finite separable morphism $f : D \to \tC_i$ such that $C \silo D^{(r)}$ for some $r \ge 1$ and $\pi_i$ factors as
\[
\pi_i : C \silo D^{(r)} \xrightarrow{F^r} D \xrightarrow{f} \tC_i \; .
\]
Write $E_i = \alpha^*_i E$. Then we have to show that $\pi^* E = \pi^*_i E_i$ is semistable of degree $0$ on $C$. By assumption, $E_i$ is strongly semistable of degree $0$ on $\tC_i$. Using Lemma \ref{t2_1} and the fact that $F$ commutes with all morphisms in characteristic $p$, we find that the pullback $f^* E_i$ under the finite, separable map $f$ is strongly semistable of degree $0$ on $D$. Hence $F^{r*} f^* E_i$ is semistable, which implies that $\pi^*_i E_i$ is semistable of degree zero.
\end{proof}

Generalizing a result by Lange and Stuhler in \cite{La-Stu}, one can show

\begin{prop}
  \label{t2_5}
If $k = \F_q$ is a finite field, then a vector bundle $E$ on the proper $k$-curve $Z$ is strongly semistable of degree zero, if and only if there exists a finite surjective morphism 
\[
\varphi : Y \to Z
\]
of proper $k$-curves such that $\varphi^* E$ is trivial. In fact, one can take $\varphi$ to be the composition
\[
\varphi : Y \xrightarrow{\Fr^s_q} Y \xrightarrow{\pi} Z
\]
of a power of the $k$-linear Frobenius morphism $\Fr_q$ (defined by the $q$-th power map on $\Oh_Y$) and a finite, \'etale and surjective morphism $\pi$.
\end{prop}

\begin{proof}
  See \cite{De-We2}, Theorem 18.
\end{proof}

Note that every vector bundle $E$ of degree zero on a smooth geometrically connected projective curve $C$ of genus $g \le 1$ is semistable. Namely, for $g = 0$, every vector bundle of degree zero is in fact trivial. For $g = 1$, the claim follows from Atiyah's classification \cite{At}: let $E = \bigoplus_i E_i$ be the decomposition of $E$ into indecomposable components. Since $E$ is semistable of degree zero, all $E_i$ have degree zero since they are subbundles and quotients. Therefore by \cite{At} Theorem 5,  we have $E_i \simeq L \otimes G$, where $L$ is a line bundle of degree zero and $G$ is an iterated extension of trivial line bundles. The pullback of $E_i$ under some Frobenius power is also of this form. Since the category of semistable vector bundles of degree $0$ on $C$ is closed under extensions and contains all line bundles of degree zero, we conclude that $E_i$ is indeed strongly semistable of degree $0$. This proves the following fact:

\begin{lemma}
  \label{t2_6}
Let $Z$ be a proper $k$-curve such that the normalizations $\tC_i$ of all irreducible components $C_i$ are geometrically connected of genus $g (\tC_i) \le 1$. Consider a vector bundle $E$ on $Z$. If all restriction $E \, |_{\tC_i}$ are semistable of degree zero, then $E$ is strongly semistable of degree zero. 
\end{lemma}

By \cite{Gie1}, for every genus $\ge 2$ there are examples of semistable vector bundles of degree zero which are not strongly semistable.

On the other hand, there are results indicating that there are ``a lot of'' strongly semistable vector bundles of degree zero. In \cite{La-Pau}, Laszlo and Pauly show that for an ordinary smooth projective curve $C$ of genus two over an algebraically closed field $k$ of characteristic two, the set of strongly semistable rank two bundles is Zariski dense in the coarse moduli space of all semistable rank two bundles with trivial determinant. See \cite{JRXY} for generalizations to higher genus.

In \cite{Du}, Ducrohet investigates the case of a supersingular smooth projective curve of genus two over an algebraically closed field $k$ with $\car (k) = 2$. It turns out that in this case all equivalence classes of semistable bundles with trivial determinant but one are in fact strongly semistable.

\section{Vector bundles on $p$-adic curves}
\label{sec:3}

Before discussing the $p$-adic case, let us recall some results in the complex case, i.e. regarding vector bundles on a compact Riemann surface $X$. Let $x \in X$ be a base point and denote by $\pi : \tX \to X$ the universal covering of $X$. Every representation $\rho : \pi_1 (X,x) \to \GL_r (\C)$ gives rise to a flat vector bundle $E_{\rho}$ on $X$, which is defined as the quotient of the trivial bundle $\tX \times \C^r$ by the $\pi_1 (X,x)$-action given by combining the natural action of $\pi_1 (X,x)$ on the first factor with the action induced by $\rho$ on the second factor. It is easily seen that every flat vector bundle on $X$ is isomorphic to some $E_{\rho}$. Regarding $E_{\rho}$ as a holomorphic bundle on $X$, a theorem of Weil \cite{W} says that a holomorphic bundle $E$ on $X$ is isomorphic to some $E_{\rho}$ (i.e. $E$ comes from a representation of $\pi_1 (X,x)$) if and only if in the decomposition $E = \bigoplus^r_{i=1} E_i$ of $E$ into indecomposable subbundles all $E_i$ have degree zero. A famous result by Narasimhan and Seshadri \cite{Na-Se} says that a holomorphic vector bundle $E$ of degree $0$ on $X$ is stable if and only if $E$ is isomorphic to $E_{\rho}$ for some irreducible unitary representation $\rho$. Hence a holomorphic vector bundle comes from a unitary representation $\rho$ if and only if it is of the form $E = \bigoplus^r_{i=1} E_i$ for stable (and hence indecomposable) subbundles of degree zero.

Now let us turn to the $p$-adic case. Let $X$ be a connected smooth projective curve over the algebraic closure $\oQ_p$ of $\Q_p$ and put $X_{\C_p} = X \otimes_{\oQ_p} \C_p$. We want to look at $p$-adic representations of the algebraic fundamental group $\pi_1 (X , x)$ where $x \in X (\C_p)$ is a base point. It is defined as follows. Denote by $F_x$ the functor from the category of finite \'etale coverings $X'$ of $X$ to the category of finite sets which maps $X'$ to the set of $\C_p$-valued points of $X'$ lying over $x$.

For $x , x' \in X (\C_p)$ we call any isomorphism $F_x \silo F_{x'}$ of fibre functors an \'etale path from $x$ to $x'$. (Note that any topological path on a Riemann surface induces naturally such an isomorphism of fibre functors.) Then the \'etale fundamental group $\pi_1 (X,x)$ is defined as
\[
\pi_1 (X,x) = \Aut (F_x) \; .
\]

The goal of our papers \cite{De-We1} and \cite{De-We2} is to associate $p$-adic representations of the \'etale fundamental group $\pi_1 (X,x)$ to certain vector bundles on $X_{\C_p}$. Let us briefly describe the main result. We call any finitely presented, proper and flat scheme $\eX$ over the integral closure $\oZ_p$ of $\Z_p$ in $\oQ_p$ with generic fibre $X$ a model of $X$. By $\eo$ we denote the ring of integers in $\C_p$, and by $k = \oF_p$ the residue field of $\oZ_p$ and $\eo$. We write $\eX_{\eo} = \eX \otimes_{\Z_p} \eo$ and $\eX_k = \eX \otimes_{\oZ_p} k$. 

\begin{defn}
  \label{t3_1}
We say that a vector bundle $E$ on $X_{\C_p}$ has strongly semistable reduction of degree zero if $E$ is isomorphic to the generic fibre of a vector bundle $\Eh$ on $\eX_{\eo}$ for some model $\eX$ of $X$, such that the special fibre $\Eh_k$ is a strongly semistable vector bundle of degree zero on the proper $k$-curve $\eX_k$.

$E$ has potentially strongly semistable reduction of degree zero if there is a finite morphism $\alpha : Y \to X$ of connected smooth projective curves over $\oQ_p$ such that $\alpha^*_{\C_p} E$ has strongly semistable reduction of degree zero on $Y_{\C_p}$.
\end{defn}

By $\eB^s$ (respectively $\eB^{ps}$) we denote the full subcategory of the category of vector bundles on $X_{\C_p}$ consisting of all $E$ with strongly semistable (respectively potentially strongly semistable) reduction of degree zero. Besides, for every divisor $D$ on $X_{\C_p}$ we define $\eB_{X_{\C_p} , D}$ to be the full subcategory of those vector bundles $E$ on $X_{\C_p}$ which can be extended to a vector bundle $\Eh$ on $\eX_{\eo}$ for some model $\eX$ of $X$, such that there exists a finitely presented proper $\oZ_p$-morphism 
\[
\pi : \Yh \longrightarrow \eX
\]
satisfying the following two properties:\\
i) The generic fibre of $\pi$ is finite and \'etale outside $D$\\
ii) The pullback $\pi^*_k \Eh_k$ of the special fibre of $\Eh$ is trivial on $\Yh_k$ (c.f. \cite{De-We2}, definition 6 and theorem 16).

Then we show in \cite{De-We2}, Theorem 17:
\[
\eB^s = \bigcup_D \eB_{X_{\C_p} , D} \; ,
\]
where $D$ runs through all divisors on $X_{\C_p}$. By \cite{De-We2}, Theorem 13, every bundle in $\eB_{X_{\C_p}, D}$ is semistable of degree zero, so that $\eB^s$ and also $\eB^{ps}$ are full subcategories of the category $\eT^{ss}$ of semistable bundles of degree zero on $X_{\C_p}$. Line bundles of degree zero lie in $\eB^{ps}$ by \cite{De-We2} Theorem 12\,a.

The main result in \cite{De-We2} is the following (c.f. \cite{De-We2}, theorem 36):

\begin{theorem}
  \label{t3_2}
Let $E$ be a bundle in $\eB^{ps}$. For every \'etale path from $x$ to $y$ in $X (\C_p)$ there is an isomorphism
\[
\rho_E (\gamma) : E_x \silo E_y
\]
of ``parallel transport'', which behaves functorially in $\gamma$. The association $E \mapsto \rho_E (\gamma)$ is compatible with tensor products, duals and internal homs of vector bundles in the obvious way. It is also compatible with $\Gal(\oQ_p / \Q_p)$-conjugation. Besides, if $\alpha : X \to X'$ is a morphism of smooth projective curves over $\oQ_p$ and $E'$ a bundle in $\eB^{ps}_{X'_{\C_p}}$, then $\rho_{\alpha^* E'} (\gamma)$ and $\rho_{E'} (\alpha_* \gamma)$ coincide, where $\alpha_* \gamma$ is the induced \'etale path on $X'$. For every $x \in X (\C_p)$ the fibre functor
\[
\eB^{ps} \longrightarrow \uVec_{\C_p} \; ,
\]
mapping $E$ to the fibre $E_x$ in the category $\uVec_{\C_p}$ of $\C_p$-vector spaces, is faithful.
\end{theorem}

In particular one obtains a continuous representation $\rho_{E,x} : \pi_1 (X,x) \to \GL_r (E_x)$. The functor $E \mapsto \rho_{E,x}$ is compatible with tensor products, duals, internal homs, pullbacks of vector bundles and $\Gal (\oQ_p / \Q_p)$-conjugation.

Let us look at two special cases of this representation: For line bundles on a curve $X$ with good reduction, $\rho$ induces a homomorphism
\[
\Pic^0_X (\C_p) \longrightarrow \Hom_{\cont} (\pi_1 (X,x) , \C^*_p)
\]
mapping $L$ to $\rho_{L,x}$. As shown in \cite{De-We1} this map coincides with the map defined by Tate in \cite{Ta} \S\,4 on an open subgroup of $\Pic^0_X (\C_p)$. Secondly, applying $\rho$ to bundles $E$ in $H^1 (X_{\C_p} , \Oh) = \Ext^1_{X_{\C_p}} (\Oh , \Oh)$, one recovers the Hodge--Tate map to $H^1 (X_{\et} , \Q_p) \otimes \C_p = \Ext^1_{\pi_1 (X,x)} (\C_p , \C_p)$, see \cite{De-We1}, corollary 8.

It follows from \cite{De-We2}, Proposition 9 and Theorem 11 that the categories $\eB^s$ and $\eB^{ps}$ are closed under tensor products, duals, internal homs and extensions. We will now prove another important property of those categories.

\begin{theorem}
  \label{t3_4}
If a vector bundle $E$ on $X_{\C_p}$ is contained in $\eB^s$ (respectively $\eB^{ps}$), then every quotient bundle of degree zero and every subbundle of degree zero of $E$ is also contained in $\eB^s$ (respectively $\eB^{ps}$).
\end{theorem}

\begin{proof}
  It suffices to show this property for the category $\eB^s$. By duality, it suffices to treat quotient bundles. So let $\tEh$ be a vector bundle with strongly semistable reduction of degree zero on $\eX_{\eo}$, where $\eX$ is a model of $X$. Denote by $E$ the generic fibre of $\tEh$, and let
  \begin{equation}
    \label{eq:1}
    0 \to E' \to E \to E'' \to 0
  \end{equation}
be an exact sequence of vector bundles $X_{\C_p}$, where $E''$ has degree zero.
By \cite{De-We2}, Theorem 5, $E'$ can be extended to a vector bundle $\Fh'$ on $\Yh_{\eo}$, where $\Yh$ is a model of $X$ such that there is a morphism $\varphi : \Yh \to \eX$ inducing an isomorphism on the generic fibres. Since $\Hom (\Fh' , \varphi^*_{\eo} \tEh) \otimes_{\eo} \C_p = \Hom (E' , E)$, we may assume that the embedding $E' \to E$ can be extended to a $\Oh_{\Yh_{\eo}}$-module homomorphism $\Fh' \to \varphi^*_{\eo} \tEh$ after changing the morphisms in the diagram (\ref{eq:1}).

Let $\Fh''$ be the quasi-coherent sheaf on $\Yh_{\eo}$ such that $\Fh' \to \varphi^*_{\eo} \tEh \to \Fh'' \to 0$ is exact. Then $\Fh''$ is of finite presentation. Note that the generic fibre of this sequence is isomorphic to the sequence (\ref{eq:1}).

Let $r$ be the rank of $E''$. The same argument as in the proof of \cite{De-We2}, Theorem 5 shows that the blowing-up
\[
\psi_{\eo} : \eZ_{\eo} \to \Yh_{\eo}
\]
of the $r$-th Fitting ideal of $\Fh''$ descends to a finitely presented morphism $\psi : \eZ \to \Yh$ inducing an isomorphism on the generic fibres. Besides, if $\Ah$ denotes the annihilator of the $r$-th Fitting ideal of $\psi^*_{\eo} \Fh''$, the sheaf $\psi^*_{\eo} \Fh'' / \Ah$ is locally free by \cite{RG} (5.4.3). Hence it gives rise to a vector bundle $\Eh''$ on $\eZ_{\eo}$ with generic fibre $E''$. Let us write $\Eh = \psi^*_{\eo} \varphi^*_{\eo} \tEh$. Then we have a natural surjective homomorphism of vector bundles $\Eh \to \Eh''$ on $\eZ_{\eo}$ extending the quotient map $E \to E''$ on the generic fibre $X_{\C_p}$.

Let $K$ be a finite extension of $\Q_p$ such that $\eZ$ descends to a proper and flat scheme $\eZ_{\eo_K}$ over the ring of integers $\eo_K$. We can choose $K$ big enough so that all irreducible components of the special fibre of $\eZ$ are defined over the residue field of $K$. Let $X_K$ be the generic fibre of $\eZ_{\eo_K}$. Then $X_K \otimes_K \oQ_p \simeq X$. The scheme $\eZ_{\eo}$ is the projective limit of all $\eZ_A = \eZ_{\eo_K} \otimes_{\eo_K} A$, where $A$ runs over the finitely generated $\eo_K$-subalgebras of $\eo$.

By \cite{EGAIV}, (8.5.2), (8.5.5), (8.5.7), (11.2.6) there exists a finitely generated $\eo_K$-subalgebra $A$ of $\eo$ with quotient field $Q \subset \C_p$ such that $\Eh \to \Eh''$ descends to a surjective homomorphism $\Eh_A \to \Eh''_A$ of vector bundles on $\eZ_A$. 

Let $x \in \spec A$ be the point corresponding to the prime ideal $A \cap \emm$ in $A$, where $\emm \subset \eo$ is the valuation ideal. If $\pi_K$ is a prime element in $\eo_K$, we have $\eo_K / (\pi_K) \subset A / A \cap \emm \subset \eo / \emm = k$, so that $A \cap \emm$ is a maximal ideal in $A$. Hence $x$ is a closed point with residue field $\kappa = \kappa (x)$ which is a finite extension of $\eo_K / (\pi_K)$ in $k$.

 By assumption, the vector bundle $\tEh_k$ on the special fibre $\eX_k$ of $\eX$ is strongly semistable of degree zero. 
By Proposition \ref{t2_4}, strong semistability is preserved under pullbacks via $k$-morphisms, so that $\Eh_k = (\varphi \verk \psi)^*_k \tEh_k$ is also strongly semistable of degree zero.

The bundle $\Eh_{\kappa} = \Eh_A \otimes_A \kappa$ satisfies $\Eh_{\kappa} \otimes_{\kappa} k \simeq \Eh_k$, hence it is strongly semistable of degree zero on $\eZ_{\kappa} = \eZ_A \otimes_A \kappa$. 

Let $C_1 , \ldots , C_r$ be the irreducible components of $\eZ_{\kappa}$ with normalizations $\tC_1 , \ldots , \tC_r$ and denote by $\alpha_i : \tC_i \to C_i \to \eZ_{\kappa}$ the natural map. Since the Euler characteristics are locally constant in the fibres of the flat and proper $A$-scheme $\eZ_A$, we find $\deg \Eh''_{\kappa} = \deg E''_Q = 0$.

By the degree formula in \cite{BLR}, 9.1, Proposition 5, $\deg (\Eh''_{\kappa})$ is a linear combination of the $\deg (\alpha^*_i \Eh''_{\kappa})$'s with positive coefficients. Since $\alpha^*_i \Eh''_{\kappa}$ is a quotient bundle of the semistable degree zero vector bundle $\alpha^*_i \Eh_{\kappa}$ on $\tC_i$, it has degree $\ge 0$. Hence for all $i$ we find $\deg (\alpha^*_i \Eh''_{\kappa}) = 0$.

Now let $\Fh$ be a vector bundle on the smooth projective curve $\tC_i$ which is a quotient of $\alpha^*_i \Eh''_{\kappa}$. Then $\Fh$ is also a quotient of the semistable degree zero bundle $\alpha^*_i \Eh_{\kappa}$, which implies $\deg (\Fh) \ge \deg \alpha^*_i \Eh_{\kappa} = 0$. This shows that $\alpha^*_i \Eh''_{\kappa}$ is semistable of degree $0$ on $\tC_i$. The same argument applies to all Frobenius pullbacks of $\alpha^*_i \Eh''_{\kappa}$, so that $\alpha^*_i \Eh''_{\kappa}$ is strongly semistable of degree $0$. Hence $\Eh''_{\kappa}$ is strongly semistable of degree zero on $\Zh_{\kappa}$. By \cite{HL}, 1.3.8 the base change $\Eh''_k = \Eh''_{\kappa} \otimes_{\kappa} k$ is also strongly semistable of degree zero. Since $\Eh''$ has generic fibre $E''$, it follows that $E''$ is indeed contained in $\eB^s$.
\end{proof}

\begin{cor}
  \label{t3_5}
$\eB^s$ and $\eB^{ps}$ are abelian categories.
\end{cor}

\begin{proof}
  Recall that $\eB^s$ and $\eB^{ps}$ are full subcategories of the abelian category $\Th^{ss}$ of semistable vector bundles of degree zero on $X_{\C_p}$. Since the trivial bundle is contained in $\eB^s$ and $\eB^{ps}$, and both categories are closed under direct sums by \cite{De-We1}, Proposition 9, they are additive.

By the theorem, $\eB^s$ and $\eB^{ps}$ are also closed under kernels and cokernels, hence they are abelian categories.
\end{proof}

\section{Tannakian categories of vector bundles}
\label{sec:4}

In this section we look at several categories of semistable vector bundles from a Tannakian point of view. Useful references in this context are \cite{De-Mi} and \cite{S} for example.

As before let $X$ be a smooth projective curve over $\oQ_p$ with a base point $x \in X (\C_p)$. We call a vector bundle on $X_{\C_p}$ polystable of degree zero if it is isomorphic to the direct sum of stable vector bundles of degree zero. Let $\Th^{ss}_{\red}$ be the strictly full subcategory of vector bundles on $X_{\C_p}$ consisting of polystable bundles of degree zero and set $\eB^{ps}_{\red} = \eB^{ps} \cap \Th^{ss}_{\red}$. Then we have the following diagram of fully faithful embeddings
\begin{equation}
  \label{eq:2}
  \begin{array}{ccc}
\eB^{ps}_{\red} & \subset & \eB^{ps} \\
\cap & & \cap \\
\Th^{ss}_{\red} & \subset & \Th^{ss}
  \end{array}
\end{equation}
Note that because of theorem \ref{t3_4}, every vector bundle $E$ in $\eB^{ps}_{\red}$ is the direct sum of stable vector bundles of degree zero contained in $\eB^{ps}$.

\begin{lemma}
  \label{t4_1}
The categories $\Th^{ss}_{\red}$ and $\eB^{ps}_{\red}$ are closed under taking subquotients in $\Th^{ss}$.
\end{lemma}

\begin{proof}
  Since $\eB^{ps}$ is closed under subquotients in $\Th^{ss}$ by theorem \ref{t3_4}, it suffices to consider $\Th^{ss}_{\red}$. A bundle $E$ in $\Th^{ss}_{\red}$ can be written as $E = \bigoplus_i E_i$ with $E_i$ stable of degree zero. Let $\varphi : E \to E''$ be a surjective map in $\Th^{ss}$. Then the images $\varphi (E_i)$ lie in $\Th^{ss}$ and the surjective map $\varphi |_{E_i} : E_i \to \varphi (E_i)$ is an isomorphism or the zero map since $E_i$ is stable. Hence we have $E'' = \sum_j E''_j$ for stable degree zero vector bundles $E''_j$ (the nonzero $\varphi (E_i)$). For every $j$ we have
\[
E''_j \cap \sum_{k \neq j} E''_k = E''_j \quad \mbox{or} \; = 0 
\]
since the bundle $E''_j$ being stable is a simple object of $\Th_{ss}$. It follows that $E''$ is the direct sum of suitably chosen $E''_j$'s and hence lies in $\Th^{ss}_{\red}$. The case of subobjects follows by duality. 
\end{proof}

Consider the fibre functor $\omega_x$ on $\Th^{ss}$ defined by $\omega_x (E) = E_x$ and $\omega_x (f) = f_x$. It induces fibre functors on the other categories as well. 

\begin{theorem}
  \label{t4_2}
{\bf a} \quad The categories $\eB^{ps}_{\red} , \eB^{ps} , \Th^{ss}_{\red}$ and $\Th^{ss}$ with the fibre functor $\omega_x$ are neutral Tannakian categories over $\C_p$.\\
{\bf b} \quad The categories $\eB^{ps}_{\red}$ and $\Th^{ss}_{\red}$ are semisimple. Every object in $\eB^{ps}$ (resp. $\Th^{ss}$) is a successive extension of objects of $\eB^{ps}_{\red}$ (resp. $\Th^{ss}_{\red}$).\\
{\bf c} \quad The natural inclusion $\eB^{ps} \subset \Th^{ss}$ is an equivalence of categories if and only if $\eB^{ps}_{\red} \subset \Th^{ss}_{\red}$ is an equivalence of categories.
\end{theorem}

\begin{proof}
  {\bf a} \quad For $\Th^{ss}$ and $\Th^{ss}_{\red}$ this is well known, see e.g. \cite{Si}, p. 29. The categories $\eB^{ps}$ and $\eB^{ps}_{\red}$ are abelian by corollary \ref{t3_5} and lemma \ref{t4_1}. It was shown in \cite{De-We2} that $\eB^{ps}$ is closed under tensor products and duals. The same follows for $\eB^{ps}_{\red} = \eB^{ps} \cap \Th^{ss}_{\red}$. Faithfulness of $\omega_x$ on $\eB^{ps}$ and $\eB^{ps}_{\red}$ follows because $\omega_x$ is faithful on $\Th^{ss}$. Alternatively a direct proof was given in \cite{De-We2} Theorem 36.\\
{\bf b} \quad Every object in $\eB^{ps}_{\red}$ and $\Th^{ss}_{\red}$ is the direct sum of simple objects since stable bundles are simple. It is well known that objects of $\Th^{ss}$ are successive extensions of stable bundles of degree zero. Since subquotients in $\Th^{ss}$ of objects in $\eB^{ps}$ lie in $\eB^{ps}$ by theorem \ref{t3_4}, the corresponding assertion for $\eB^{ps}$ follows.\\
{\bf c} \quad This is a consequence of {\bf b} because both $\Th^{ss}$ and $\eB^{ps}$ are closed under extensions, c.f. \cite{De-We2}. 
\end{proof}

Let
\begin{equation}
  \label{eq:3}
  \vcenter{\xymatrix{
G^{ps}_{\red} & G^{ps} \ar[l] \\
G^{ss}_{\red} \ar[u] & G^{ss} \ar[u] \ar[l]
}}
\end{equation}
be the diagram of affine group schemes over $\C_p$ corresponding to diagram (\ref{eq:2}) by Tannakian duality.

\begin{prop}
  \label{t4_3}
All morphisms in (\ref{eq:3}) are faithfully flat. The connected components of $G^{ps}_{\red}$ and $G^{ss}_{\red}$ are pro-reductive.
\end{prop}

\begin{proof}
  The following is known \cite{De-Mi} Proposition 2.21: \\
A fully faithful $\otimes$-functor $F : \Ch \to \Dh$ of neutral Tannakian categories over a field $k$ of characteristic zero induces a faithfully flat morphism $F^* : G_{\Dh} \to G_{\Ch}$ of the Tannakian duals if and only if we have: Every subobject in $\Dh$ of an object $F (C)$ for some $C$ in $\Ch$ is isomorphic to $F (C')$ for a subobject $C'$ of $C$. \\
This criterion can be verified immediately for the functors in (\ref{eq:2}) by using either theorem \ref{t3_4} or lemma \ref{t4_1}.
The second assertion of the proposition follows from \cite{De-Mi} Proposition 2.23 and Remark 2.28. 
\end{proof}

Let $\Th_{\fin}$ be the category of vector bundles on $X_{\C_p}$ which are trivialized by a finite \'etale covering of $X_{\C_p}$. Noting that $\pi_1 (X_{\C_p} , x) \silo \pi_1 (X,x)$ is an isomorphism it follows that $\Th_{\fin}$ is equivalent to the category of representations of $\pi_1 (X,x)$ with open kernels on finite dimensional $\C_p$-vector spaces $V$, c.f. \cite{La-Stu}, 1.2. Such a representation factors over a finite quotient $G$ of $\pi_1 (X,x)$ and the corresponding bundle in $\Th_{\fin}$ is $E = X'_{\C_p} \times^G \bV$. Here $\alpha : X' \to X$ is the Galois covering corresponding to the quotient $\pi_1 (X,x) \to G$ and $\bV$ is the affine space over $\C_p$ corresponding to $V$. With the fibre functor $\omega_x$ the category $\Th_{\fin}$ is neutral Tannakian over $\C_p$ with Tannaka dual
\[
\pi_1 (X,x)_{/ \C_p} = \varprojlim_N (\pi_1 (X,x) / N)_{/ \C_p} \; .
\]
Here $N$ runs over the open normal subgroups of $\pi_1 (X,x)$ and for a finite (abstract) group $H$ we denote by $H_{/ \C_p}$ the corresponding constant group scheme. Using Maschke's theorem it follows that $\Th_{\fin}$ is semisimple.

\begin{prop}
  \label{t4_4}
The category $\Th_{\fin}$ is a full subcategory of $\eB^{ps}_{\red}$. The induced morphism $G^{ps}_{\red} \twoheadrightarrow \pi_1 (X,x)_{/ \C_p}$ is faithfully flat.
\end{prop}

\begin{proof}
  Using \cite{De-We2} Proposition 37, one sees that $\Th_{\fin}$ is a full subcategory of $\eB^{ps}$. If $V$ is a finite dimensional $\C_p$-vector space on which $\pi_1 (X,x)$ acts via a finite quotient, we can decompose $V$ into a direct sum $V = \bigoplus_i V_i$ of irreducible representations using Maschke's theorem. Using an algebraic isomorphism of $\C_p$ with $\C$ and Narasimhan--Seshadri theory it follows that the degree zero vector bundle $E_i = X'_{\C_p} \times^G \bV_i$ is stable. Hence $E = \bigoplus E_i$ is an object of $\Th^{ss}_{\red}$ and hence of $\eB^{ps}_{\red} = \eB^{ps} \cap \Th^{ss}_{\red}$. The next assertion follows from fully faithfullness of $\Th_{\fin} \hookrightarrow \eB^{ps}_{\red}$ since $\eB^{ps}_{\red}$ is semisimple, c.f. \cite{De-Mi}, Remark 2.29.
\end{proof}

Consider a Galois covering $\alpha : X' \to X$ with group $\Gal (X' / X)$ of smooth projective curves over $\oQ_p$ and choose a point $x' \in X' (\C_p)$ above $x \in X (\C_p)$. Let us write $\Ch_X$ for any of the categories $\eB^{ps}_{\red} , \eB^{ps} , \Th^{ss}_{\red}$ and $\Th^{ss}$ of vector bundles on $X_{\C_p}$. The pullback functor $\alpha^* : \Ch_X \to \Ch_{X'}$ is a morphism of neutral Tannakian categories over $\C_p$ commuting with the fibre functors $\omega_x$ and $\omega_{x'}$. Let $i : G_{X'} \to G_X$ be the morphism of Tannaka duals induced by $\alpha^*$. We also need the faithfully flat homomorphism obtained by composition:
\[
q : G_X \twoheadrightarrow \pi_1 (X , x)_{/ \C_p} \twoheadrightarrow \Gal (X' / X)_{/ \C_p} \; .
\]
Here the second arrow is determined by our choice of $x'$. Note that every $\sigma$ in $\Gal (X' / X)$ induces an automorphism $\sigma^*$ of $\Ch_{X'}$ and hence an automorphism $\sigma : G_{X'} \to G_{X'}$ of group schemes over $\C_p$. 

\begin{lemma}
  \label{t4_5} 
There is a natural exact sequence of affine group schemes over $\C_p$
\[
1 \to G_{X'} \xrightarrow{i} G_X \xrightarrow{q} \Gal (X' / X)_{/ \C_p} \to 1 \; .
\]
\end{lemma}

\begin{proof}
  Every bundle $E'$ in $\Ch_{X'}$ is isomorphic to a subquotient of $\alpha^* (E)$ for some bundle $E$ in $\Ch_X$. Namely, thinking of $E'$ as a locally free sheaf, the sheaf $E = \alpha_* E'$ is locally free again and we have $\alpha^* E \cong \bigoplus_{\sigma} \sigma^* E'$ where $\sigma$ runs over $\Gal (X' / X)$. Incidentally, $E$ lies in $\Ch_X$ because $\alpha^* E \cong \bigoplus_{\sigma} \sigma^* E'$ lies in $\Ch_{X'}$. This is clear for $\Ch = \Th^{ss}$ or $\eB^{ps}$. For $\Th^{ss}_{\red}$ and hence for $\eB^{ps}_{\red}$ it follows from \cite{HL} Lemma 3.2.3. 

It follows from \cite{De-Mi} Proposition 2.21 (b) that $i$ is a closed immersion. By descent the category $\Ch_X$ is equivalent to the category of bundles in $\Ch_{X'}$ equipped with a $\Gal (X'/X)$-operation covering the one on $X'$. In other words, the category of representations of $G_X$ is equivalent to the category of representations of $G_{X'}$ together with a $\Gal (X' / X)$-action, i.e. a transitive system of isomorphisms $\sigma^* \rho = \rho \verk \sigma \to \rho$ for all $\sigma$ in $\Gal (X' / X)$. Hence, $G_X$ is an extension of $G_{X'}$ by $\Gal (X' / X)_{/ \C_p}$ inducing the above $\Gal (X' / X)$-action on $G_{X'}$. (Because such an extension has the same $\otimes$-category of representations.) In particular, the sequence in the lemma is exact.
\end{proof}

For a commutative diagram of Galois coverings
\[
\xymatrix{
X'' \ar[rr] \ar[dr] & & X' \ar[dl]\\
 & X
}
\]
and the choice of points $x' \in X' (\C_p)$ and $x'' \in X'' (\C_p)$ over $x$ we get a commutative diagram of affine group schemes over $\C_p$:
\[
\xymatrix{
1 \ar[r] & G_{X''} \ar[r] \ar@{_{(}->}[d] & G_X \ar[r] \ar@{=}[d] & \Gal (X'' / X)_{/ \C_p} \ar[r] \ar@{>>}[d] & 0 \\
1 \ar[r] & G_{X'} \ar[r] & G_X \ar[r] & \Gal (X' / X)_{/ \C_p} \ar[r] & 0 \; .
}
\]
Passing to the projective limit, we get an exact sequence
\begin{equation}
  \label{eq:4}
  1 \longrightarrow \varprojlim_{X'} G_{X'} \longrightarrow G_X \longrightarrow \pi_1 (X , x)_{/\C_p} \longrightarrow 1 \; .
\end{equation}
Right exactness follows from propositions \ref{t4_3} and \ref{t4_4}.

\begin{theorem}
  \label{t4_6}
We have a commutative diagram
\[
\xymatrix{
1 \ar[r] & (G^{ss}_{\red})^0 \ar[r] \ar@{>>}[d] & G^{ss}_{\red} \ar[r] \ar@{>>}[d] & \pi_1 (X,x)_{/ \C_p} \ar@{=}[d] \ar[r] & 1 \\
1 \ar[r] & (G^{ps}_{\red})^0 \ar[r] & G^{ps}_{\red} \ar[r] & \pi_1 (X,x)_{/ \C_p} \ar[r] & 1 
}
\]
In particular $\pi_1 (X,x)_{/ \C_p}$ is the common group scheme of connected components of both $G^{ss}_{\red}$ and $G^{ps}_{\red}$. Moreover we have:
\[
(G^{ss}_{\red})^0 = \varprojlim_{X'} G^{ss}_{\red, X'} \quad \mbox{and} \quad (G^{ps}_{\red})^0 = \varprojlim_{X'} G^{ps}_{\red, X'} \; .
\]
Here $X' / X$ runs over a cofinal system of pointed Galois covers of $(X, x)$.
\end{theorem}

\begin{proof}
  Let $\Ch_X$ denote either $\eB^{ps}_{\red , X}$ or $\Th^{ss}_{\red, X}$ and let $G_X$ be its Tannaka dual. The exact sequence (\ref{eq:4}) implies that $G^0_X \subset \varprojlim_{X'} G_{X'}$. Hence it suffices to show that $\varprojlim_{X'} G_{X'}$ is connected. The category of finite dimensional representations of $\varprojlim_{X'} G_{X'}$ on $\C_p$-vector spaces is $\varinjlim_{X'} \Ch_{X'}$. In order to show that $\varprojlim_{X'} G_{X'}$ is connected, by \cite{De-Mi} Corollary 2.22 we have to prove the following:

{\bf Claim} Let $A$ be an object of $\varinjlim_{X'} \Ch_{X'}$. Then the strictly full subcategory $[[A]]$ of $\varinjlim_{X'} \Ch_{X'}$ whose objects are isomorphic to subquotients of $A^N , N \ge 0$ is not stable under $\otimes$ unless $A$ is isomorphic to a trivial bundle.

\begin{proofob}
  Let $[[A]]$ be stable under $\otimes$. The category $\varinjlim_{X'} \Ch_{X'}$ is semisimple since $(\varprojlim_{X'} G_{X'})^0 = \varprojlim_{X'} G^0_{X'}$ is pro-reductive. Hence we may decompose $A$ into simple objects $A = A_1 \oplus \ldots \oplus A_s$. By assumption, for every $j \ge 1$ the object $A^{\otimes j}_1$ is isomorphic to a subquotient of $A^N$ for some $N = N (j)$. The same argument as in the proof of Lemma \ref{t4_1} shows that up to isomorphism the subquotients of $NA := A^N$ have the form $m_1 A_1 \oplus \ldots \oplus m_s A_s$ for integers $m_i \ge 0$. Hence we get isomorphisms where $\sum$ means ``direct sum'':
\[
A^{\otimes j}_1 \cong \sum^s_{i=1} m_{ij} A_i \quad \mbox{for} \; 1 \le j \le r \; .
\]
Here $M = (m_{ij})$ is an $s \times r$-matrix over $\Z$. Fixing some $r > s$ there is a relation with integers $c_j$, not all zero:
\[
\sum^r_{j=1} c_j (m_{1j} , \ldots , m_{sj})^t = 0 \; .
\]
This gives the relation
\[
\sum^r_{j=1} c^+_j (m_{1j} , \ldots , m_{sj})^t = \sum^r_{j=1} c^-_j (m_{1j} , \ldots , m_{sj})^t
\]
where $c^+_j = \max \{ c_j , 0 \}$ and $c^-_j = - \min \{ c_j , 0 \}$. ``Left multiplication'' with $(A_1 , \ldots , A_s)$ gives isomorphisms 
\[
\sum^r_{j=1} c^+_j \sum^s_{i=1} m_{ij} A_i \cong \sum^r_{j=1} c^-_j \sum^s_{i=1} m_{ij} A_i \; ,
\]
and hence 
\[
\sum^r_{j=1} c^+_j A^{\otimes j}_1 \cong \sum^r_{j=1} c^-_j A^{\otimes j}_1 \; .
\]
For the polynomials $P^{\pm} (T) = \sum^r_{j=1} c^{\pm}_j T^j$ with coefficients in $\Z^{\ge 0}$ we have $P^+ \neq P^-$ and:
\[
P^+ (A_1) \cong P^- (A_1) \; .
\]
Let $E_1$ be a bundle in $\Ch_{X'}$ representing $A_1$ in $\varinjlim_{X'} \Ch_{X'}$. Then we have an isomorphism 
\[
P^+ (\beta^* E_1) \cong P^- (\beta^* E_1) 
\]
of vector bundles on a suitable Galois cover $\beta : X'' \to X'$. A theorem of Weil, c.f. \cite{W} Ch. III or \cite{N}, now implies that $\beta^* E_1$ and hence $E_1$ is trivialized by a finite \'etale covering of $X'$. Hence $A_1$, the class of $E_1$ is isomorphic in $\varinjlim_{X'} \Ch_{X'}$ to a trivial bundle. The same argument applies to $A_2 , A_3 , \ldots $ Hence $A$ is isomorphic to a trivial bundle as well. This proves the claim and hence the theorem.
\end{proofob}
\end{proof}

We now determine the structure of $G^{\abb}$ for $G = G^{ps}_{\red}$ and $G = G^{ss}_{\red}$. This group is pro-reductive and abelian, hence diagonalizable and therefore determined by its character group
\[
X (G^{\abb}) = \Mor_{\C_p} (G^{\abb} , \Ge_m) \; .
\]
The characters of $G^{\abb}$ correspond to isomorphism classes of one-dimensional representations of $G$ i.e. to isomorphism classes of degree zero line bundles in $\eB^{ps}_{\red}$ resp. $\Th^{ss}_{\red}$. Since both categories contain all degree zero line bundles we get
\[
X (G^{\abb}) = \Pic^0_X (\C_p)
\]
and hence
\begin{eqnarray*}
  G^{\abb} & = & \Hom (\Pic^0_X (\C_p) , \Ge_{m , \C_p}) \\
& = & \varprojlim_A \Hom (A , \Ge_{m,\C_p}) \; .
\end{eqnarray*}
Here $A$ runs over the finitely generated subgroups of $\Pic^0_X (\C_p)$. A similar argument using the fact that $(G^{ps}_{\red})^0$ resp. $(G^{ss}_{\red})^0$ is the Tannaka dual of $\varinjlim_{X'} \eB^{ps}_{\red , X'}$ resp. $\varinjlim_{X'} \Th^{ss}_{\red, X'}$ shows the following: For $G$ as above, the group $(G^0)^{\abb}$ is diagonalizable with character group
\[
X ((G^0)^{\abb}) = \varinjlim_{X'} \Pic^0_{X'} (\C_p) \; .
\]
Note that the right hand group is torsionfree because line bundles of finite order become trivial in suitable finite \'etale coverings. This corresponds to the fact that $(G^0)^{\abb}$ is connected. We can therefore write as well:
\[
X ((G^0)^{\abb}) = \varinjlim_{X'} (\Pic^0_{X'} (\C_p) / \tors)
\]
and hence
\[
(G^0)^{\abb} = \varprojlim_{X'} \Hom (\Pic^0_{X'} (\C_p) / \tors , \Ge_{m,\C_p}) \; .
\]
This is a pro-torus over $\C_p$.

In particular we have seen that $(G^{ps}_{\red})^0$ and $(G^{ss}_{\red})^0$ have the same maximal abelian quotient. Incidentally we may compare $(G,G^0)$ with $(G^0,G^0)$:\\
General theory gives an exact sequence:
\[
1 \longrightarrow (G,G^0) / (G^0,G^0) \longrightarrow (G^0)^{\abb} \longrightarrow (G^{\abb})^0 \longrightarrow 1 \; .
\]
The sequence of character groups is
\[
1 \longrightarrow \Pic^0_X (\C_p) / \tors \longrightarrow \varinjlim_{X'} (\Pic^0_{X'} (\C_p) / \tors) \longrightarrow X ((G , G^0) / (G^0,G^0)) \longrightarrow 1 \; .
\]
If $X$ is an elliptic curve, it follows that the first arrow is an isomorphism, so that $(G, G^0) = (G^0,G^0)$. If the genus satisfies $g (X) \ge 2$ then by the Riemann Hurwitz formula, the middle group is infinite dimensional and in particular $(G , G^0) / (G^0 , G^0)$ contains a non-trivial pro-torus.

We proceed with some remarks on the structure of $G^0$ which follow from the general theory of reductive groups. Let $C$ be the neutral component of the center of $G^0$. Then we have
\[
G^0 = C \cdot (G^0 , G^0) \; .
\]
Here $C$ is a pro-torus and $(G^0 , G^0)$ is pro-semisimple. The projection $C \to (G^0)^{\abb}$ is faithfully flat and its kernel is a commutative pro-finite groupscheme $H$. In the exact sequence:
\[
1 \longrightarrow X ((G^0)^{\abb}) \longrightarrow X (C) \longrightarrow X (H) \longrightarrow 1
\]
the group $X ((G^0)^{\abb})$ is divisible because $\Pic^0_{X'} (\C_p)$ is divisible. Since $X ((G^0)^{\abb})$ is also torsionfree it is a $\Q$-vector space. The group $X (H)$ being torsion it follows that $X (C) \otimes \Q = X ((G^0)^{\abb})$ canonically and $X (C) \cong X ((G^0)^{\abb}) \oplus X (H)$ non-canonically. It would be interesting to determine $X (H)$ for both $G = G^{ps}_{\red}$ and $G^{ss}_{\red}$. 

We end with a remark on the Tannaka dual $G_E$ of the Tannaka subcategory of $\eB^{ps}_{\red}$ generated by a vector bundle $E$ in $\eB^{ps}_{\red}$. The group $G_E$ is a subgroup of $\bGL_{E_x}$ the linear group over $\C_p$ of the $\C_p$-vector space $E_x$. It can be characterized as follows: The group $G_E (\C_p)$ consists of all $g$ in $\GL (E_x)$ with $g (s_x) = s_x$ for all $n,m \ge 0$ and all sections $s$ in
\[
\Gamma (X_{\C_p} , (E^*)^{\otimes n} \otimes E^{\otimes m}) = \Hom_{X_{\C_p}} (E^{\otimes n} , E^{\otimes m}) \; .
\]
Here $g (s_x)$ means the extension of $g$ to an automorphism $g$ of $(E^*_x)^{\otimes n} \otimes E^{\otimes m}_x$ applied to $s_x$.

Consider the representation attached to $E$ by theorem \ref{t3_2}:
\[
\rho_{E,x} : \pi_1 (X,x) \longrightarrow \GL (E_x) \; .
\]
Its image is contained in $G_E (\C_p)$ because the functor $F \mapsto \rho_{F,x}$ on $\eB^{ps}_{\red}$ is compatible with tensor products and duals and maps the trivial line bundle to the trivial representation. Hence $G_E$ contains the Zariski closure of $\Imm \rho_{E,x}$ in $\bGL_{E_x}$. It follows from a result by Faltings \cite{Fa} that the faithful functor $F \mapsto \rho_{F,x}$ is in fact fully faithful. If $\rho_{E,x}$ is also semisimple then a standard argument shows that $G_E$ is actually equal to the Zariski closure of $\Imm \rho_{E,x}$.
\begin{minipage}{6cm}
Mathematisches Institut\\
Einsteinstr. 62\\
48149 M\"unster\\
Germany\\
deninger@math.uni-muenster.de
\end{minipage}
\hspace*{\fill} \begin{minipage}{8cm}
Fachbereich Mathematik\\
Pfaffenwaldring 57\\
70569 Stuttgart\\
Germany\\
werner@mathematik.uni-stuttgart.de
\end{minipage}
\end{document}